\documentclass[smallextended,referee,envcountsect]{svjour3}
\smartqed
\usepackage{graphicx}
\usepackage{epstopdf}
\usepackage{cite}
\usepackage[misc]{ifsym}
\usepackage[breaklinks,
           pdfstartview=FitH,
           colorlinks=true,
           citecolor=blue]{hyperref} 
\usepackage{graphicx}
\usepackage{latexsym,amsmath,amssymb}
\journalname{}

\begin{document}

\title{Solvability of a Regular Polynomial Vector Optimization Problem without Convexity\thanks{This work was partially supported by the National Natural Science  Foundation of China (11471230).}}

\author{Danyang Liu$^1$\and Rong Hu$^2$ \and Yaping Fang$^1$}

\institute{Danyang Liu \\
              394898525@qq.com  \\
           Rong Hu\\
              ronghumath@aliyun.com\\
           Yaping Fang,  Corresponding author \\
              ypfang@scu.edu.cn\\
              $^1$College of Mathematics, Sichuan University, Chengdu, Sichuan, P.R. China\\
              $^2$College of Applied Mathematics, Chengdu University of Information Technology, Chengdu, Sichuan, P.R. China
}


\maketitle

\begin{abstract}
In this paper we consider the solvability of  a non-convex regular polynomial vector optimization problem on  a nonempty closed set. We introduce  regularity conditions for the polynomial vector optimization problem and study  properties and characterizations of the regularity conditions. Under the regularity conditions, we study  nonemptiness and boundedness of the solution sets of the problem. As a consequence, we establish two Frank-Wolfe type theorems for the non-convex polynomial vector optimization problem. Finally, we investigate the solution stability of the non-convex regular polynomial vector optimization problem.
\end{abstract}

\section{Introduction}

Throughout,  $\mathbf{R}^{n}$ denotes the $n$-dimensional Euclidean space with  the norm $\|\cdot\|$ and the inner  $\langle \cdot,\cdot\rangle$,  and $\mathbf{R}^{n}_{+}:=\{x=(x_1,\cdots,x_n)\in \mathbf{R}^{n}:x_i \ge 0, i=1, \cdots, n\}$. In this paper we are interested in the following polynomial vector optimization problem on  $K$:
$$ \text{PVOP}(K, f): \qquad\mathbf{R}^{s}_{+}-Min_{x\in K}f(x),$$
where $f=(f_{1},  \dots, f_{s}): \mathbf{R}^{n} \mapsto \mathbf{R}^{s}$ is  a vector polynomial such that each component function  $f_{i}$ is a polynomial with its degree $ \mbox{ deg } f_i=d_{i}$, and $K\subseteq \mathbf{R}^{n}$ is a nonempty closed set (not necessarily convex set or semi-algebraic set  \cite{HHV,BR}).

\par Recall that a point $x^{*}\in K$ is a \emph{Pareto efficient solution} of $\text{PVOP}(K, f)$ if $$f(x)-f(x^{*})\notin - \mathbf{R}^{s}_{+}\backslash \{0\}, \quad \forall x\in K$$
and   $x^{*}\in K$ is a \emph{weak Pareto efficient solution} of PVOP$(K, f)$ if $$f(x)-f(x^{*})\notin - \mbox{ int }\mathbf{R}^{s}_{+}, \quad \forall x\in K.$$
The Pareto efficient solution set and the weak Pareto efficient solution set of $\text{PVOP}(K, f)$ are denoted by $SOL^{s}(K, f)$ and $SOL^{w}(K, f)$ respectively. Clearly, $SOL^{s}(K, f)\subseteq SOL^{w}(K, f)$.  When $s=1$,  $\text{PVOP}(K, f)$ collapses to the  polynomial scalar optimization problem:
$$\mbox{PSOP}(K, f):\qquad \emph{Min}_{x\in K}f(x),$$
whose  solution set is denoted by $SOL(K, f)$.

For the polynomial scalar optimization problem, regularity condition has been  used  in \cite{AQ} to investigate the nonemptiness of the  solution set and the continuity of the solution mapping for a quadratic programming problem.  Hieu \cite{HV1}  established a Frank-Wolfe type theorem for a polynomial scalar optimization problem on a nonempty closed set by proving the  nonemptiness and boundedness of its solution set when the objective function is bounded from below on the constraint set   and the regularity condition holds.  Hieu et al. \cite{HV2} proved that the solution set of an optimization problem corresponding to a polynomial complementarity problem is nonempty and compact  by using the regularity condition of the polynomial complementarity problem. Kim et al.  \cite{DTN} proved the existence of Pareto efficient solutions of an unconstrained polynomial vector optimization problem when the  Palais-Smale-type condition holds and the image of the objective function has a bounded section. When  $K$ is a convex semi-algebraic  set and  $f$ is convex, Jiao et al. \cite{LGJ} proved that  $\text{PVOP}(K, f)$ has a Pareto efficient solution if and only if  the image $f(K)$ of $f$ has a nonempty  bounded section.  Inspired by the above works, in this paper, we study  the nonemptiness and boundedness of the solution sets of $\text{PVOP}(K, f)$ without assuming any convexity of the objective function.

The rest of this paper is organized as follows: In Section 2, we present some notations and preliminary results.  In Section 3, we  give some sufficient and necessary conditions for regularity conditions. In Section 4, we  discuss local properties of the regularity conditions. Section 5 is devoted to the study of solvability of $\text{PVOP}(K, f)$ under regularity conditions.  In Section 6, we discuss the solution stability under the regularity conditions. Finally, we makes a concluding remark in Section 7.

\section{Preliminaries}

In this section, we recall some concepts and results that will be used in this paper. A nonempty subset $C$ of $\mathbf{R}^{n}$ is called a cone if $tx\in C$ for any $x\in C$ and any $t>0$. Given a nonempty closed set $K\subset  \mathbf{R}^{n}$, the \emph{asymptotic cone} $K_{\infty}$  of $K$ is defined by
$$K_{\infty}=\{v\in \mathbf{R}^{n}: \mbox{there}\ \mbox{exist}\ t_{k}\rightarrow +\infty\ \mbox{and}\ x_{k}\in K\ \mbox{such}\ \mbox{that}\ \lim_{k\rightarrow +\infty}\frac{x_{k}}{t_{k}}=v\}.$$
As known, $K_{\infty}$ is a closed cone and $(K_{\infty})_{\infty}=K_{\infty}$, and $K$ is bounded if and only if $K_{\infty}=\{0\}$. If $K$ is a convex set, then $K_{\infty}$ is also a convex cone and $K_{\infty}=\mbox{Rec}\ K$, where $\mbox{Rec}\ K$ is the \emph{recession cone} of $K$ defined by $$\mbox{Rec}\ K=\{v\in \mathbf{R}^{n}: x+tv\in K ,\, \forall x\in K, \forall t\geq 0\}.$$
It is known that $K+K_{\infty}=K$. The above results can been found in  \cite{RC,AA}.

\begin{definition}
Let $p=(p_{1}, \dots, p_{s}): \mathbf{R}^{n} \mapsto \mathbf{R}^{s}$ be a vector polynomial  with $  \mbox{ deg }p_{i}=d_{i}$, $i=1,  \dots, s$. We say that  $p^{\infty}$ is the \emph{vector recession polynomial} (or the \emph{vector leading term}) of $p$, where
$$p^{\infty}(x)=(p^{\infty}_{1}(x), p^{\infty}_{2}(x), \dots, p^{\infty}_{s}(x)) \emph{ and } p^{\infty}_{i}(x)=\lim_{\lambda\to +\infty}\frac{p_{i}(\lambda x)}{\lambda^{d_{i}}},\quad\forall x\in \mathbf{R}^{n}.$$

\end{definition}

\begin{remark}\label{recpoly}
	When $s=1$, $p^{\infty}$ is a recession polynomial of $p$ (see \cite{HV1}).
\end{remark}

\begin{definition}(See e.g. \cite{BT1,BT2,DTN})
A function $h: \mathbf{R}^{n} \mapsto \mathbf{R}^{s}$ is said to be bounded from below on $K$ if  there exists $r\in \mathbf{R}^{s}$ such that $h(x)\in r+\mathbf{R}^{s}_{+}$ for all $x\in K$. Clearly, $h$ is  bounded from below on $K$ if and only if its component $h_i$ is  bounded from below on $K$ for all $i$.
\end{definition}

\begin{definition} We say that $\text{ PVOP }(K, f)$ is \emph{weakly regular} (resp.  \emph{strongly  regular}) if  $SOL^{s}(K_{\infty}, f^{\infty})$ (resp. $SOL^{w}(K_{\infty}, f^{\infty})$) is bounded.
\end{definition}
\begin{remark}
Clearly,  strong regularity implies  weak regularity.  When $s=1$, both weak  regularity  and strong  regularity  coincide with the regularity  in \cite[Definition 2.1]{HV1}.
\end{remark}

 The following result plays an important role in establishing the existence of the Pareto efficient solutions  of $\text{PVOP}(K,F)$.

\begin{lemma}\cite[Proposition 13]{MG} \label{scalars}
	Given  $\lambda\in \mathrm{ int }\ \mathbf{R}^{s}_{+}$ and $x_{0}\in K$, define  $g(x)=\sum^{s}_{j=1}\lambda_{j}f_{j}(x)$ and $G_{x_{0}}=\{x\in  K: f_{i}(x)\leq f_{i}(x_{0}), i=1, 2, \dots, s\}$. If $x^{*}\in SOL(G_{x_{0}}, g)$, then $x^{*}\in SOL^{s}(K, f)$.
\end{lemma}

In what follows we always assume the each component polynomial $f_i$ of the objective function $f$ has a degree $d_{i}\geq 1$.

\section{Regularity of $\text{PVOP}(K, f)$}

In this section, we shall discuss properties and characterizations of  regularity  of $\text{PVOP}(K, f)$.

\subsection{Conditions for regularity}

In this subsection we shall show that the  regularity  of $\text{PVOP}(K, f)$ is closely related to the regularity  of $\text{PSOP}(K, f_i)$.  To do so, we first give a characterization of $SOL^{w}(K_{\infty}, f^{\infty})=\emptyset$.

\begin{proposition}\label{necessary}
	$SOL^{w}(K_{\infty}, f^{\infty})=\emptyset$  if and only if $0\notin SOL^{w}(K_{\infty}, f^{\infty})$.
\end{proposition}
{\it Proof}   We only need to prove the sufficiency.	 Suppose that $0\notin SOL^{w}(K_{\infty}, f^{\infty})$. Then there exists $v_{0}\in K_{\infty}\backslash\{0\}$ such that
$$f^{\infty}(v_{0})-f^{\infty}(0)=f^{\infty}(v_{0})\in - \mbox{ int }\mathbf{R}^{s}_{+},$$
which yields
$$f^{\infty}_{i}(v_{0})<0,\quad i=1,  \dots, s.$$
Let $v\in K_{\infty}$. It follows that
$$f^{\infty}_{i}(tv_{0})-f^{\infty}_{i}(v)=t^{d_{i}}f^{\infty}_{i}(v_{0})-f^{\infty}_{i}(v)<0, \quad i=1,\cdots, s$$
for all sufficiently large $t>0$.  Since $v\in K_{\infty}$ is arbitrary,  $SOL^{w}(K_{\infty}, f^{\infty})=\emptyset$.
\qed

\begin{example}\label{emptyex}
	Consider the vector polynomial $f=(f_{1}, f_{2})$  with
$$f_{1}(x_{1}, x_{2})=x^{3}_{1}-x^{2}_{1}x_{2}-3x_{1}+2x_{2}+1, f_{2}(x_{1}, x_{2})=-x^{2}_{2}-x_{1}x_{2}+x_{1}-1$$
and
$$K=\{(x_{1}, x_{2})\in \mathbf{R}^{2}: x_{1}\geq 0, x_{2}-x_{1}\geq 0\}.$$
It is easy to verify that $K=K_{\infty}$,  $f^{\infty}_{1}(x_{1}, x_{2})=x^{3}_{1}-x^{2}_{1}x_{2}$, and $f^{\infty}_{2}(x_{1}, x_{2})=-x^{2}_{2}-x_{1}x_{2}$.  Then $0\notin SOL^{w}(K_{\infty}, f^{\infty})$ since $f^{\infty}_{1}(1, 2)=-1<0=f^{\infty}_{1}(0, 0)$ and $f^{\infty}_{2}(1, 2)=-6<0=f^{\infty}_{2}(0, 0)$. By Proposition \ref{necessary}, $SOL^{w}(K_{\infty}, f^{\infty})=\emptyset$.
\end{example}

 \begin{proposition}\label{cones}
 	If $SOL^{w}(K_{\infty}, f^{\infty})\neq\emptyset$, then $SOL^{w}(K_{\infty}, f^{\infty})$ is a nonempty cone.
 \end{proposition}
 {\it Proof}
 	 Since $SOL^{w}(K_{\infty}, f^{\infty})\neq\emptyset$, by Proposition \ref{necessary}, $0\in SOL^{w}(K_{\infty}, f^{\infty})$. Suppose on the contrary that  there exist $v_{0}\in SOL^{w}(K_{\infty}, f^{\infty})$ and $t_{0}>0$ such that $t_{0}v_{0}\notin SOL^{w}(K_{\infty}, f^{\infty})$. Then  there exists $v_{1}\in K_{\infty}$ such that
	    \begin{equation}\label{eq11}
	    	f^{\infty}_{i}(v_{1})-f^{\infty}_{i}(t_{0}v_{0})<0, \quad  i=1,  \dots, s.
	    \end{equation}
Dividing the  both sides of the inequality (\ref{eq11}) by $t^{d_{i}}_{0}$, we obtain $f^{\infty}_{i}(\frac{v_{1}}{t_{0}})-f^{\infty}_{i}(v_{0})<0$ for all $i\in \{1, \dots, s\}$. This reaches a contradiction to  $v_{0}\in SOL^{w}(K_{\infty}, f^{\infty})$ since $\frac{v_{1}}{t_{0}}\in K_{\infty}$. Thus, $SOL^{w}(K_{\infty}, f^{\infty})$ is a nonempty cone.
\qed

\begin{remark}\label{RE}
	By Proposition \ref{cones}, \text{PVOP}$(K, f)$ is strongly  regular if and only if $SOL^{w}(K_{\infty}, f^{\infty})$ is empty or $SOL^{w}(K_{\infty}, f^{\infty})=\{0\}$.
\end{remark}

\begin{proposition}\label{notbounded}
	If $SOL^{w}(K_{\infty}, f^{\infty})=\emptyset$, then $f_{i}$ is unbounded from below on $K$ for all $i\in \{1, \dots, s\}$.
\end{proposition}
{\it Proof}
Suppose on the contrary that  there exists $i_{0}\in \{1, 2, \dots, s\}$ such that $f_{i_{0}}$ is bounded from below on $K$. Then  $SOL(K_{\infty}, f^{\infty}_{i_{0}})=\emptyset$ since $SOL(K_{\infty}, f^{\infty}_{i_{0}})\subseteq SOL^{w}(K_{\infty}, f^{\infty})$.
By Lemma \ref{necessary}, there exists $v_{0}\in K_{\infty}\backslash\{0\}$ such that $f^{\infty}_{i_{0}}(v_{0})<f^{\infty}_{i_{0}}(0)=0$.  Since $v_{0}\in K_{\infty}\backslash\{0\}$, there exist $t_k>0$ with $t_{k}\rightarrow +\infty$ and $x_{k}\in K$ such that $t^{-1}_{k}x_{k}\rightarrow v_{0}$ as $k\rightarrow +\infty$. Since $f_{i_{0}}$ is bounded from below on $K$, there exists a constant $l$ such that
$$\frac{f_{i_{0}}(x_{k})}{t^{d_{i_{0}}}_{k}}\geq \frac{l}{t^{d_{i_{0}}}_{k}}.$$
Letting $k\rightarrow +\infty$, we obtain $f^{\infty}_{i_{0}}(v_{0})\geq 0$ which is a contradiction to $f^{\infty}_{i_{0}}(v_{0})<0$.
\qed

When $f_i$ is bounded from below on $K$ for some $i\in\{1, \cdots, s\}$, by Proposition \ref{notbounded} and  Proposition \ref{necessary},  we know that $0\in SOL^{w}(K_{\infty}, f^{\infty})$.  In the following proposition, we further show  $0\in SOL^{s}(K_{\infty}, f^{\infty})$ when each component $f_i$ is bounded from below on $K$.

\begin{proposition}\label{weakpro}
	If $f$ is bounded from below on $K$, then $0\in SOL^{s}(K_{\infty}, f^{\infty})$.
\end{proposition}
{\it Proof}
	Suppose that $f$ is bounded from below on $K$. Then $f_i$ is  bounded from below on $K$ for all $i$. By Proposition \ref{notbounded} and  Proposition \ref{necessary}, $f^{\infty}_{i}(x)\geq f^{\infty}_{i}(0)=0$ for all $x\in K_{\infty}$ and  for all $i\in \{1,  \dots, s\}$. This means $0\in SOL^{s}(K_{\infty}, f^{\infty})$.
\qed

The following result shows that the  weak (strong) regularity of $\text{ PVOP }(K,f)$ is closely related to the regularity of $\text{ PVOP } (K, f_i)$.

\begin{theorem}\label{gpro}
	\item(i) $SOL^{w}(K_{\infty}, f^{\infty})=\{0\}$ if and only if $SOL(K_{\infty}, f^{\infty}_{i})= \{0\}$ for all $i\in \{1, \dots, s\}$.
	\item(ii) If $SOL^{s}(K_{\infty}, f^{\infty})=\{0\}$, then $SOL(K_{\infty}, f^{\infty}_{i})\neq\emptyset$ for all $i\in \{1, \dots, s\}$.
\end{theorem}
{\it Proof}
\emph{(i)}  Assume that $SOL^{w}(K_{\infty}, f^{\infty})=\{0\}$.  Since
$$SOL(K_{\infty}, f^{\infty}_{i})\subseteq SOL^{w}(K_{\infty}, f^{\infty}),$$ it suffices to show $SOL(K_{\infty}, f^{\infty}_{i})\neq\emptyset$ for all $i\in \{1, \dots, s\}$. Suppose on the contrary that there exists $i_{0}\in \{1, \dots, s\}$ such that $SOL(K_{\infty}, f^{\infty}_{i_{0}})=\emptyset$ which is equivalent to $0\notin SOL(K_{\infty}, f^{\infty}_{i_{0}})$ (by Lemma \ref{necessary}).  Then there exists $v_{0}\in K_{\infty}\backslash \{0\}$ such that $f^{\infty}_{i_{0}}(v_{0})<0$.  Let $\lambda\in \mbox{ int }\mathbf{R}^{s}_{+}$. Consider the function
\begin{equation}\label{18-fyp1}
F_{\lambda}(x)=\sum^{s}_{i=1}\lambda_{i}f^{\infty}_{i}(x)\; \text{ and }\;
S_{v_{0}}=\{x\in K_{\infty}: f^{\infty}(x)\leq f^{\infty}(v_{0})\}.
 \end{equation}
Then $S_{v_{0}}$ is nonempty and closed. Next we shall show that $S_{v_{0}}$ is bounded. If not, then there exists the sequence $\{x_{k}\}\subseteq S_{v_{0}}$ such that $\|x_{k}\|\to +\infty$ as $k\to +\infty$. Without loss of generality, we assume that $\|x_{k}\|\neq 0$ and $\frac{x_{k}}{\|x_{k}\|}\to x^{*}\in K_{\infty}\backslash \{0\}$. Since   $x_{k}\in  S_{v_{0}}$, we have
 $$f^{\infty}_{i}(x_{k})\leq f^{\infty}_{i}(v_{0}),\quad\forall i=1,  \dots, s.$$
Dividing the both sides of the above inequality by $\|x_{k}\|^{d_{i}}$ and letting $k\rightarrow +\infty$, we get
$$
	f^{\infty}_{i}(x^{*})\leq 0=f^{\infty}_{i}(0), \quad\forall i=1,  \dots, s,
$$
which together with  $0\in SOL^{w}(K_{\infty}, f^{\infty})$ yields $x^{*}\in SOL^{w}(K_{\infty}, f^{\infty})$, a contradiction. Thus, $S_{v_{0}}$ is bounded.  By the known  \emph{Weierstrass'  theorem}, we have $SOL(S_{v_{0}}, F_{\lambda})\neq\emptyset$. It follows from  Lemma \ref{scalars} that
$$\emptyset\ne SOL(S_{v_{0}}, F_{\lambda})\subseteq SOL^{s}(K_{\infty}, f^{\infty})\subseteq SOL^{w}(K_{\infty}, f^{\infty})=\{0\},$$
which implies  $SOL(S_{v_{0}}, F_{\lambda})=\{0\}$, and so $0\in S_{v_{0}}$. By the definition of $ S_{v_{0}}$, we get $f^{\infty}_{i_{0}}(v_{0})\geq 0$, a contradiction to $f^{\infty}_{i_{0}}(v_{0})<0$. Therefore, $SOL(K_{\infty}, f^{\infty}_{i})\neq\emptyset$ for all  $i\in \{1,  \dots, s\}$.

For the converse, assume that $SOL(K_{\infty}, f^{\infty}_{i})= \{0\}$ for all $i\in \{1,  \dots, s\}$ and there exists $v_{0}\in SOL^{w}(K_{\infty}, f^{\infty})\backslash\{0\}$. Then
$$f^{\infty}(0)-f^{\infty}(v_{0})\notin - \mbox{ int }\mathbf{R}^{s}_{+},$$
which implies that  $f^{\infty}_{i_0}(v_{0})\le f^{\infty}_{i_0}(0)$ for some  $i_0\in \{1,  \dots, s\}$. This reaches  a contradiction to  $SOL(K_{\infty}, f^{\infty}_{i_0})= \{0\}$.

 \emph{(ii)} Suppose on the contrary that there exists $i_{0}\in \{1, \dots, s\}$ such that $SOL(K_{\infty}, f^{\infty}_{i_{0}})=\emptyset$. By Lemma \ref{necessary}, $0\notin SOL(K_{\infty}, f^{\infty}_{i_{0}})$.
Then there exists $v_{0}\in K_{\infty}\backslash \{0\}$ such that $f^{\infty}_{i_{0}}(v_{0})<0$.  Let $\lambda\in \mbox{ int }\mathbf{R}^{s}_{+}$. By similar arguments as in the proof of (i), we have
$$\emptyset\ne SOL(S_{v_{0}}, F_{\lambda})\subseteq SOL^{s}(K_{\infty}, f^{\infty})=\{0\},$$
where $S_{v_{0}}$ and $F_{\lambda}$ are defined as in \eqref{18-fyp1}. The rest is same as the one of (i), and so we omit it.
\qed

\begin{remark}
The following example shows that $SOL^{s}(K_{\infty}, f^{\infty})=\{0\}$ does not imply  $SOL(K_{\infty}, f^{\infty}_{i})= \{0\}$ for each $i\in \{1,  \dots, s\}$.

\end{remark}

\begin{example}
	Consider the vector polynomial $f=(f_{1}, f_{2})$ with
	$$f_{1}(x_{1}, x_{2})=x^{3}_{1}-x_{2}+1, f_{2}(x_{1}, x_{2})=x^{3}_{2}-x_{1}-1$$
	and the constraint set
	$$K=\{(x_{1}, x_{2})\in \mathbf{R}^{2}: x_{1}\geq 0, x_{2}\geq 0\}.$$
It is easy to verify that $K=K_{\infty}$, $f^{\infty}_{1}(x_{1}, x_{2})=x^{3}_{1}$ and $f^{\infty}_{2}(x_{1}, x_{2})=x^{3}_{2}$. Clearly, $SOL^{s}(K_{\infty}, f^{\infty})=\{0\}$. However, $SOL(K_{\infty}, f^{\infty}_{1})= \{(x_{1}, x_{2})\in \mathbf{R}^{2}: x_{1}=0, x_{2}\geq 0\}\neq\{0\}$ and $SOL(K_{\infty}, f^{\infty}_{2})= \{(x_{1}, x_{2})\in \mathbf{R}^{2}: x_{1}\geq 0, x_{2}=0\}\neq\{0\}$.
\end{example}

\begin{proposition}\label{eqcondition}
$\mathrm{PVOP}(K, f)$ is strongly  regular and $f$ is bounded from below on $K$ if and only if  $SOL(K_{\infty}, f^{\infty}_{i})= \{0\}$ for all $i\in \{1, \dots, s\}$.
\end{proposition}

{\it Proof}
The sufficiency follows immediately from Proposition \ref{weakpro} and  Theorem \ref{gpro} \emph{(i)}.

Next we prove the sufficiency. By  Theorem \ref{gpro}\emph{(i)}, $\text{ PVOP }(K, f)$ is strongly  regular. Suppose on the contrary that there exists $i_0\in \{1, \dots, s\}$ such that $f_{i_0}$ is not bounded from below on $K$. Let $x\in K$. Then there exists $\{x_{k}\}^{\infty}_{k=1}\subseteq K$ such that
\begin{equation}\label{fyp1230a}
f_{i_0}(x_{k})\leq -k\leq f_{i_0}(x)
\end{equation}
for all sufficiently large $k$. We claim that $\{x_{k}\}^{\infty}_{k=1}\subseteq K$ is unbounded. Indeed, if not, then we may assume that  $\|x_{k}\|\to +\infty$ and $\frac{x_{k}}{\|x_{k}\|}\to v\in K_{\infty}\backslash \{0\}$  as $k\to +\infty$. Dividing the both sides  of (\ref{fyp1230a}) by $\|x_{k}\|^{d_{i}}$ and letting $k\rightarrow +\infty$, we get $f^{\infty}_{i}(v)\leq 0$, a contradiction to $SOL(K_{\infty}, f^{\infty}_{i_0})=\{0\}$. Hence, $\{x_{k}\}^{\infty}_{k=1}$ is bounded. Without loss of generality, we may assume that $x_{k}\to x^{*}\in K$ as $k\to +\infty$.  It follows from (\ref{fyp1230a})  that
$$f_{i_0}(x^{*})=\lim_{k\to +\infty}f_{i_0}(x_{k})\leq \lim_{k\to +\infty}-k=-\infty,$$
 a contradiction.  Thus, $f$ is bounded from below on $K$.
\qed

\subsection{Regularity and  $R_{0}$-property}

It has been shown in \cite{Oett,FYP,FangH,Gowda,Lop} that  $R_0$-property plays an important role in studying the compactness of the solution sets of complementarity problems as well as  the upper continuity of the solution mappings.  In this subsection, we shall show that the weak (strong) regularity of $\mbox{ PVOP }(K,f)$ is closely related to the $R_0$-property of the following  weak vector complementarity problem \cite{CGY,Yangxq}
$$\mbox{WVCP}(K_{\infty}, \nabla f^{\infty}): \ \mbox{Find}\ x\in K_{\infty} \, \mbox{such that}\, \nabla f^{\infty}(x)\in (K_{\infty})^{w+}_{\mathbf{R}^{s}_+}\,\mbox{and}\, \langle \nabla f^{\infty}(x), x\rangle\notin \mbox{ int } \mathbf{R}^{s}_+,$$
where $\nabla f^{\infty}=(\nabla f^{\infty}_{1}, \nabla f^{\infty}_{2}, \dots, \nabla f^{\infty}_{s})$ is the gradient of $f^{\infty}$,
$$(K_{\infty})^{w+}_{\mathbf{R}^{s}_+}=\{v\in L(\mathbf{R}^{n}, \mathbf{R}^{s}): \langle v, x\rangle\notin -\mbox{ int }\mathbf{R}^{s}_+, \, \forall x\in K_{\infty}\}$$
is the weak dual cone of $K_{\infty}$ with respect to $\mathbf{R}^{s}_+$, and $L(\mathbf{R}^{n}, \mathbf{R}^{s})$ is the space of all bounded linear operators from $\mathbf{R}^{n}$ to $\mathbf{R}^{s}$.  We denote the solution set of $\mbox{WVCP}(K_{\infty}, \nabla f^{\infty})$ by  $SOL_{WVCP}(K_{\infty},\nabla f^{\infty})$. Recall that $\mbox{WVCP}(K_{\infty}, \nabla f^{\infty})$ is of \emph{ type $R_{0}$} \cite{FYP} if $SOL_{WVCP}(K_{\infty},\nabla f^{\infty})=\{0\}$.

\begin{theorem}\label{sufficiently}
	Assume that $K$ is convex.  If $SOL^{w}(K_{\infty}, f^{\infty})=\{0\}$, then $\mbox{WVCP}(K_{\infty}, \nabla f^{\infty})$ is of \emph{ type $R_{0}$}. Conversely,  if $\mbox{WVCP}(K_{\infty}, \nabla f^{\infty})$ is of \emph{ type $R_{0}$}, then $\mbox{ PVOP }(K, f)$ is strongly regular.
	\end{theorem}
{\it Proof}
It is clear that $0\in SOL_{WVCP}(K_{\infty},\nabla f^{\infty})$.

Suppose that $SOL^{w}(K_{\infty}, f^{\infty})=\{0\}$. By Theorem \ref{gpro}\emph{(i)}, we have $$SOL(K_{\infty}, f_i^{\infty})=\{0\},  i=1,\cdots, s.$$
 This implies $f^{\infty}_{i}(v)>0$ for all $v\in K_{\infty}\backslash\{0\}, i=1,  \dots, s$. This together with the Euler's Homogeneous Function Theorem yields
$$\langle\nabla f^{\infty}(v), v\rangle=(\langle\nabla f^{\infty}_{1}(v),v\rangle, \dots, \langle\nabla f^{\infty}_{s}(v), v\rangle)=(d_{1}f^{\infty}_{1}(v),  \dots, d_{s}f^{\infty}_{s}(v))\in \mbox{ int }\mathbf{R}^{s}_{+}$$
for all $v\in K_{\infty}\backslash\{0\}$, where $d_i$ is the degree of $f_i, i=1,\cdots, s$.
As a consequence, $SOL_{WVCP}(K_{\infty},\nabla f^{\infty})=\{0\}$ and so $\text{WVCP}(K_{\infty}, \nabla f^{\infty})$ is of \emph{ type $R_{0}$}.

For the converse, suppose $\mbox{WVCP}(K_{\infty}, \nabla f^{\infty})$ is of \emph{ type $R_{0}$}. To complete the proof, it suffices to show $SOL^{w}(K_{\infty}, f^{\infty})=\{0\}$ (by Remark \ref{RE}).  Suppose  on the contrary that there exists $v^*\in SOL^{w}(K_{\infty}, f^{\infty})\backslash\{0\}$. By \cite[Theorem 4]{LKY}, we get
$$\langle\nabla f^{\infty}(v^*), x-v^*\rangle\notin -\mbox{ int }\mathbf{R}^{s}_{+},\quad\forall x\in K_{\infty}.$$
This implies  $v\in SOL_{WVCP}(K_{\infty},\nabla f^{\infty})$ since $K_{\infty}$ is a closed convex cone.  This reaches a contradiction.
\qed

The following example illustrates the conclusion of Theorem \ref{sufficiently}.

\begin{example}
	Consider the vector polynomial $f=(f_{1}, f_{2})$ with $f_{1}(x_{1}, x_{2})=x^{2}_{1}+x^{2}_{2}, f_{2}(x_{1}, x_{2})=x^{2}_{2}$ and
$$K=\{(x_{1}, x_{2})\in \mathbf{R}^{2}: x_{2}\geq x_{1}\geq 0\}.$$
Then $K=K_{\infty}$, $\nabla f^{\infty}_{1}(x_{1}, x_{2})=(2x_{1}, 2x_{2})$, and $\nabla f^{\infty}_{2}(x_{1}, x_{2})=(0, 2x_{2})$. It is easy to verify that $0\in SOL^{w}(K_{\infty}, f^{\infty})$. Let $x=(x_{1}, x_{2})\in K_{\infty}$ be such that $\langle\nabla f^{\infty}(x), x\rangle\notin \mbox{ int }\mathbf{R}^{2}_{+}$. It follows that $(2x^{2}_{1}+2x^{2}_{2}, 2x^{2}_{2})\notin \mbox{ int }\mathbf{R}^{2}_{+}$ This implies $x=(0, 0)$.  As a consequence,  $\mbox{WVCP}(K_{\infty},\nabla f^{\infty})$ is of type $R_{0}$. By Theorem \ref{sufficiently},  $\text{PVOP}(K, f)$ is strongly regular.
\end{example}

\begin{remark}\label{3.10}
 Assume that $K$ is convex and $f$ is bounded from below on $K$. Then by Theorem \ref{sufficiently} and Proposition \ref{weakpro}, $SOL_{WVCP}(\nabla f^{\infty}, K_{\infty})=\{0\}\Rightarrow SOL^{s}(K_{\infty}, f^{\infty})=\{0\}$. The following example shows that the converse  is not true in general.
\end{remark}

\begin{example}\label{3.7}
	Consider the vector polynomial $f=(f_{1}, f_{2})$ with $f_{1}(x_{1}, x_{2})=x^{2}_{1}, f_{2}(x_{1}, x_{2})=x^{2}_{2}$ and
$$K=\{(x_{1}, x_{2})\in \mathbf{R}^{2}: x_{2}\geq 0, x_{1}\geq 0\}.$$
Clearly, $K$ is convex, $K=K_{\infty}$, and $f$ is bounded from below on $K$. It is easy to verify that  $f_{1}(x_{1}, x_{2})=f^{\infty}_{1}(x_{1}, x_{2}), f_{2}(x_{1}, x_{2})=f^{\infty}_{2}(x_{1}, x_{2})$,  $\nabla f^{\infty}_{1}(x_{1}, x_{2})=(2x_{1}, 0)$, $\nabla f^{\infty}_{2}(x_{1}, x_{2})=(0, 2x_{2})$, and $SOL^{s}(K_{\infty}, f^{\infty})=\{0\}$. Let $x_{0}=(0, 1)$. Then $\langle\nabla f^{\infty}(x_{0}), x_{0}\rangle=(0, 2)\notin \mbox{ int }\mathbf{R}^{s}_{+}$ and $\langle\nabla f^{\infty}(x_{0}), w\rangle=(0, 2w_{2})\notin -\mbox{int}\mathbf{R}^{2}_{+}$ for all $w=(w_1,w_2)\in K_{\infty}$. This means $x_{0}\in SOL_{WVCP}(K_{\infty},\nabla f^{\infty})\backslash\{0\}$, and so $\text{WVCP}(K_{\infty},\nabla f^{\infty})$ is not of type $R_{0}$.
\end{example}

\section{Local Properties of Regularity Conditions}

In this section, we investigate  local properties of (strong) weak regularity of $\text{ PVOP }(K, f)$. Given an integer $d$, in what follows, we always  let $\mathbf{P}_{d}$ denote the family of all polynomials of degree at most $d$, and
	$$X^{n}_{d}(x):=(1, x_{1}, \dots, x_{n}, x^{2}_{1}, \dots, x^{2}_{n}, \dots, x^{d}_{1}, x^{d-1}_{1}x_{2}, x^{d-1}_{1}x_{3}, \dots, x^{d}_{2}, x^{d-1}_{2}x_{3}, x^{d-1}_{2}x_{3}, \dots, x^{d}_{n}),$$
whose components are listed by  the lexicographic ordering. The dimension of $\mathbf{P}_{d}$ is  denoted by $\kappa_d$. Then, for each polynomial $p\in \mathbf{P} _{d}$,  there exists a unique $\alpha\in \mathbf{R}^{\kappa_d}$ such that $p(x)=\langle \alpha, X^{n}_{d}(x)\rangle$. $\mathbf{P} _{d}$  can be endowed with a norm $\|p\|:=\|\alpha\|=\sqrt{\alpha^{2}_{1}+\dots+\alpha^{2}_{\kappa_d}}$. Let $p^k\in \mathbf{P} _{d}$ with $p^k\to p\in \mathbf{P} _{d}$ and $x^k\in \mathbf{R}^n$ with $x^k\to x\in \mathbf{R}^n$. It is easy to verify that $(p^k)^{\infty}\to p^{\infty}$ and $p^k(x^k)\to p(x)$ as $k\to +\infty$.

Given  $\mathbf{d}=(d_{1},  \dots, d_{s})\in \mathbf{R}^{s}$ with $d_i$ being an integer, $i=1,\cdots, s$,  let $\mathbf{P}_{\mathbf d}=\mathbf{P}_{d_{1}}\times\dots\times\mathbf{P}_{d_{s}}$. Denoted by $\mathbf{GR}^{\mathbf d}_{w}$  (resp. $\mathbf{GR}^{\mathbf d}_{s}$) the family of  all vector polynomials $p$  with $\text{ deg }p_{i}=d_{i}, i=1, \dots, s$, such that $\text{ PVOP }(K, p)$ is  strongly  (resp. weakly)  regular. Then we have the following results.

\begin{proposition}
	$\mathbf{GR}^{\mathbf{d}}_{s}$ and $\mathbf{GR}^{\mathbf{d}}_{w}$ are nonempty.
\end{proposition}
{\it Proof}
	We only need to prove that $\mathbf{GR}^{\mathbf{d}}_{s}$ is nonempty since $\mathbf{GR}^{\mathbf{d}}_{s}\subseteq \mathbf{GR}^{\mathbf{d}}_{w}$. If $K$ is bounded, then $K_{\infty}=\{0\}$. In this case $SOL^{w}(K_{\infty}, f^{\infty})=\{0\}$, and so PVOP$(K, f)$ is strongly  regular  for every vector polynomial $f$. Suppose that $K$ is unbounded. Then there exists $x^{*}=(x^{*}_{1},  \dots, x^{*}_{n})\in K_{\infty}\backslash \{0\}$. Without loss of generality, we suppose that  $x^{*}_{i_{0}}\neq 0$. Consider the vector polynomial $f=(f_{1},  \dots, f_{s}): \mathbf{R}^{n} \mapsto \mathbf{R}^{s}$ with $f_{i}(x)=-(x^{*}_{i_{0}}x_{i_{0}})^{d_{i}},  i=1,\cdots, s$. Then $f_{i}(x)$ is a polynomials $f$ of degree $d_i$ and  $f_{i}(tx^{*})=-(x^{*}_{i_{0}})^{2d_{i}}t^{d_{i}}\to -\infty$ as $t\to +\infty$.  As a consequence,   $SOL^{w}(K_{\infty}, f^{\infty})=\emptyset$, and so $f\in \mathbf{GR}^{\mathbf{d}}_{s}$.
\qed

\begin{proposition}\label{open}
	 $\mathbf{GR}^{\mathbf{d}}_{s}$ is open in $\mathbf{P}_{\mathbf{d}}$.
\end{proposition}
{\it Proof}
	We shall prove that $\mathbf{P}_{\mathbf{d}}\backslash\mathbf{GR}^{\mathbf{d}}_{s}$ is  closed in $\mathbf{P}_{\mathbf{d}}$. Let $\{f^{k}\}\subseteq \mathbf{P}_{\mathbf{d}}\backslash\mathbf{GR}^{\mathbf{d}}_{s}$ with $f^{k}=(f^{k}_{1},  \dots, f^{k}_{s})$ such that $f^{k}=(f^{k}_{1},  \dots, f^{k}_{s})\rightarrow f=(f_{1},  \dots, f_{s})$ as $k\to +\infty$.  We can suppose that   $\mbox {deg }f_{i}=d_{i}$ for all $i\in \{1, 2, \dots, s\}$ since  $f\notin \mathbf{GR}^{\mathbf{d}}_{s}$ when  $\mbox{ deg }f_{i_{0}}<d_{i_{0}}$ for some  $i_{0}\in \{1, \dots, s\}$, where $d_i$  is the $i$-th component of $\mathbf{d}$.
Since $SOL^{w}(K_{\infty}, (f_{k})^{\infty})$ is unbounded for all $k$, there exists  $x_{k}\in SOL^{w}(K_{\infty}, (f_{k})^{\infty})$ such that $\|x_k\|\to +\infty$.  Without loss of generality, we assume that  $\frac{x_{k}}{\|x_{k}\|}\to x^{*}\in K_{\infty}\backslash \{0\}$.  We claim that   $x^{*}\in SOL^{w}(K_{\infty}, f^{\infty})$. Indeed, if not, then there exists $v\in K_{\infty}$ such that
\begin{equation}\label{fypa1}
 f^{\infty}_i(v)< f^{\infty}_i(x^*),\quad i=1,\cdots, s.
\end{equation}
Since  $x_{k}\in SOL^{w}(K_{\infty}, (f_{k})^{\infty})$, we have
$$(f_{k})^{\infty}(\|x_{k}\|v)-(f_{k})^{\infty}(x_{k})\notin -\mbox{ int }\mathbf{R}^{s}_{+}.$$
Then for each $k$, there exists $i_{k, v}\in \{1,  \dots, s\}$ such that
	$$(f^{k}_{i_{k, v}})^{\infty}(\|x_{k}\|v)-(f^{k}_{i_{k, v}})^{\infty}(x_{k})\geq 0.$$
Since the set $\{1, 2, \dots, s\}$ is finite, without loss of generality, we suppose that there exists $i_{0, v}\in \{1, \dots, s\}$ such that
$$(f^{k}_{i_{0, v}})^{\infty}(\|x_{k}\|v)-(f^{k}_{i_{0, v}})^{\infty}(x_{k})\geq 0,\quad\forall k.$$
Since   $(f^{k}_{i})^{\infty}\to f^{\infty}_{i}$ as $k\to +\infty$, dividing the both sides of the above inequality by $\|x_{k}\|^{d_{i_{0, v}}}$ and letting $k\to +\infty$, we get
$$f^{\infty}_{i_{0, v}}(v)\geq f^{\infty}_{i_{0, v}}(x^{*}).$$
This reaches a contradiction to $(\ref{fypa1})$, and so $x^{*}\in SOL^{w}(K_{\infty}, f^{\infty})\backslash\{0\}$. By Proposition \ref{cones}, $SOL^{w}(K_{\infty}, f^{\infty})$ is unbounded.
	Hence, $\mathbf{P}_{\mathbf{d}}\backslash\mathbf{GR}^{\mathbf{d}}_{s}$ is closed.
\qed

\begin{remark}\label{questions}
	When $s=1$, Proposition \ref{open} reduces to  \cite[Lemma 4.1]{HV1}.
\end{remark}

In the following result, we shall show that strong regularity of a polynomial vector optimization problem  remains stable under a small perturbation.
\begin{theorem}\label{local regu}
The following conclusions hold:
	\item{(i)} If $SOL^{w}(K_{\infty}, f^{\infty})=\{0\}$, then there exists $\epsilon>0$ such that $SOL^{w}(K_{\infty}, g^{\infty})=\{0\}$ for all $g\in \mathbf{P}_{\mathbf{d}}$ satisfying $\|g-f\|<\epsilon$;
	\item{(ii)} If $SOL^{w}(K_{\infty}, f^{\infty})=\emptyset$, then there exists $\epsilon>0$ such that $SOL^{w}(K_{\infty}, g^{\infty})=\emptyset$ for all $g\in \mathbf{P}_{\mathbf{d}}$ satisfying $\|g-f\|<\epsilon$.
\end{theorem}
 {\it Proof}
Since $\mathbf{GR}^{\mathbf{d}}_{s}$ is open in $\mathbf{P}_{\mathbf{d}}$  (by Proposition \ref{open}) and $f\in\mathbf{GR}^{\mathbf{d}}_{s}$, there exists an open ball $\mathbf{B}(f, \delta)\subseteq\mathbf{GR}^{\mathbf{d}}_{s}$ such that either $SOL(K_{\infty}, g^{\infty})=\{0\}$ or $SOL(K_{\infty}, g^{\infty})=\emptyset$ for all $g\in \mathbf{B}(f, \delta)$.

\emph{(i)} It suffices to show that there exists $\epsilon\in (0,\delta)$ such that $SOL(K_{\infty}, g^{\infty})=\{0\}$ for all $g\in \mathbf{B}(f, \epsilon)$ when $SOL^{w}(K_{\infty}, f^{\infty})=\{0\}$. Suppose on the contrary that  for any $\epsilon\in (0,\delta)$, there exists $g^{\epsilon}\in \mathbf{P}_{\mathbf{d}}$ with $\|g^{\epsilon}-f\|<\epsilon$  such that $SOL^{w}(K_{\infty}, (g^{\epsilon})^{\infty})=\emptyset$.
By Lemma \ref{necessary},  there exists $x_{\epsilon}\in K_{\infty}\backslash\{0\}$ such that
 	\begin{equation}\label{3.4}
 		(g_i^{\epsilon})^{\infty}(x_{\epsilon})<(g_i^{\epsilon})^{\infty}(0)=0,\quad i=1, \cdots, s.
 	\end{equation}
Since $g^{\epsilon}\to f$ as $\epsilon\to 0$, we have $(g^{\epsilon})^{\infty}\to f^{\infty}$ as $\epsilon\to 0$. Without loss of generality, we assume that $\frac{x_{\epsilon}}{\|x_{\epsilon}\|}\to x^*\in K_{\infty}\backslash \{0\}$ as $\epsilon\to 0$. Since $g^{\epsilon}\in \mathbf{B}(f, \epsilon)\subset  \mathbf{GR}^{\mathbf{d}}_{s}$, we get $\mbox{ deg }(g_i^{\epsilon})^{\infty}=d_i$. Dividing the both sides of (\ref{3.4}) by $\|x_{\epsilon}\|^{d_{i}}$ and letting $\epsilon\rightarrow 0$, we get
$$f_{i}^{\infty}(x^*)\leq 0, \quad i=1,\cdots, s,$$
which  reaches a contradiction to $SOL(K_{\infty}, f^{\infty})= \{0\}$.

\emph{(ii)}  It suffices to show that there exists $\epsilon\in (0,\delta)$ such that $SOL(K_{\infty}, g^{\infty})=\emptyset$ for all $g\in \mathbf{B}(f, \epsilon)$ when $SOL^{w}(K_{\infty}, f^{\infty})=\emptyset$. Suppose on the contrary that  for any $\epsilon\in (0,\delta)$, there exists $g^{\epsilon}\in \mathbf{P}_{\mathbf{d}}$ with $\|g^{\epsilon}-f\|<\epsilon$  such that $SOL^{w}(K_{\infty}, (g^{\epsilon})^{\infty})=\{0\}$.
By \emph{(i)} of Theorem \ref{gpro}, we have
 	\begin{equation*}\label{3.5}
 		0=(g_i^{\epsilon})^{\infty}(0)\leq (g_i^{\epsilon})^{\infty}(x),\quad i=1, \cdots, s
 	\end{equation*}
for any $x\in K_{\infty}$. Since $g^{\epsilon}\to f$ as $\epsilon\to 0$, we have $(g^{\epsilon})^{\infty}\to f^{\infty}$ as $\epsilon\to 0$. Letting $\epsilon\rightarrow 0$ in the above inequality, we get
$$f_{i}^{\infty}(0)\leq f_{i}^{\infty}(x), \quad i=1,\cdots, s.$$
Since $x\in K_{\infty}$ is arbitrary, again from \emph{(i)} of Theorem \ref{gpro} we get $0\in SOL(K_{\infty}, f^{\infty})$, a contradiction.
 \qed

Observe that $f^{\infty}=(f+g)^{\infty}$ for all $g=(g_{1},  \dots, g_{s})\in \mathbf{P}_{\mathbf{d}}$ with $\text{ deg }g_{i}< \text{ deg }f_{i}$, $i=1, \dots, s$. As a consequence, we have the following result.

\begin{proposition}\label{perturb} Let $f=(f_1,\cdots, f_s): \mathbf{R}^{n} \mapsto \mathbf{R}^{s}$ be a  vector polynomial. Then for any vector polynomial $g=(g_{1}, \dots, g_{s})$ with $\mathrm{ deg } \ g_{i}<\mathrm{ deg } \ f_{i}, i=1, \cdots, s$, the following conclusions hold:	
\begin{itemize}
\item[(i)]  If $SOL^{s}(K_{\infty}, f^{\infty})=\{0\}$, then $SOL^{s}(K_{\infty}, (f+g)^{\infty})= \{0\}$.

\item[(ii)]  If $SOL^{s}(K_{\infty}, f^{\infty})=\emptyset$, then $SOL^{s}(K_{\infty}, (f+g)^{\infty})= \emptyset$.

\item[(iii)] If $SOL^{w}(K_{\infty}, f^{\infty})=\{0\}$, then $SOL^{w}(K_{\infty}, (f+g)^{\infty})= \{0\}$.

\item[(iv)] If $SOL^{w}(K_{\infty}, f^{\infty})=\emptyset$, then $SOL^{w}(K_{\infty}, (f+g)^{\infty})= \emptyset$.
\end{itemize}
\end{proposition}

The following result is a direct consequence of Theorem \ref{local regu} and Proposition \ref{perturb}.

\begin{corollary}
Let $f=(f_1,\cdots, f_s): \mathbf{R}^{n} \mapsto \mathbf{R}^{s}$ be a vector polynomial. Then for any vector polynomial $g=(g_{1}, \dots, g_{s})$ with $\mathrm{ deg }\  g_{i}<\mathrm{ deg }\ f_{i}, i=1,\cdots, s$, the following conclusions hold:	

\begin{itemize}
\item[(i)]  If $\mathrm{PVOP}(K, f)$ is weakly regular, then $\mathrm{PVOP}(K, f+g)$ is weakly  regular.	
\item[(ii)] If $\mathrm{PVOP}(K, f)$ is strongly regular, then $\mathrm{PVOP}(K, f+g)$ is strongly regular.
\end{itemize}
\end{corollary}

\section{Existence Results for $\text{PVOP} (K, f)$ with  Regularity}

In this section we shall study emptiness and boundedness of the solution sets of $\text{PVOP} (K, f)$ under the regularity condition.

\subsection{The weak regularity case}
 First, we give a necessary condition for the existence of Pareto efficient solutions of  $\text{PVOP}(K, f)$.

\begin{theorem}\label{exist0} Assume that one of the following conditions hold:
\begin{itemize}
\item[(i)] $SOL^{s}(K_{\infty}, f^{\infty})=\{0\}$.
\item[(ii)] there exists $i_{0}\in \{1, \dots, s\}$ such that $SOL(K_{\infty}, f^{\infty}_{i_{0}})=\{0\}$.
\end{itemize}
Then $SOL^{s}(K, f)$ is nonempty.
\end{theorem}
{\it Proof}
Let $\lambda\in \mbox{ int }\mathbf{R}^{s}_{+}$ and   $x_{0}\in K$. Define $g_{\lambda}(x)=\sum^{s}_{i=1}\lambda_{i}f_{i}(x)$ and  $G_{x_{0}}=\{x\in K: f(x)\leq f(x_{0})\}$. Clearly, $G_{x_{0}}$ is nonempty and closed. We assert that $G_{x_{0}}$ is bounded. If not, then there exists  $\{x_{k}\}\subset G_{x_{0}}$ such that $\|x_{k}\|\to +\infty$  and  $\frac{x_{k}}{\|x_{k}\|}\to \bar x\in K_{\infty}\backslash \{0\}$ as $k\to +\infty$. It follows from $x_{k}\in G_{x_{0}}$  that
$$f_{i}(x_{k})\leq f_{i}(x_{0}),\quad i=1,\cdots, s.$$
Dividing the both sides of the above inequality  by $\|x_{k}\|^{d_{i}}$ and letting $k\rightarrow +\infty$, we get
	$$f^{\infty}_{i}(\bar x)\leq f^{\infty}_{i}(0)=0, \quad i=1,\cdots, s,$$
a contradiction to  $SOL^{s}(K_{\infty}, f^{\infty})=\{0\}$ as well as $SOL(K_{\infty}, f_{i_0}^{\infty})=\{0\}$.  Hence, $G_{x_{0}}$ is compact.  By the $\emph{Weierstrass}'$ Theorem, $SOL(G_{x_{0}}, g_{\lambda})\neq\emptyset$. Since $SOL(G_{x_{0}}, g_{\lambda})\subseteq SOL^{s}(K, f)$ (by Lemma \ref{scalars}), we have $SOL^{s}(K, f)\neq\emptyset$.
\qed

The following example shows that   $SOL^{s}(K, f)$ may  be unbounded when $SOL^{s}(K_{\infty}, f^{\infty})=\{0\}$.

\begin{example}
		Consider the vector polynomial $f=(f_{1}, f_{2})$ with $f_{1}(x_{1}, x_{2})=x^{2}_{1}-x_{2}-1, f_{2}(x_{1}, x_{2})=x^{3}_{2}+1$ and
$$K=\{(x_{1}, x_{2})\in \mathbf{R}^{2}: x_{2}\geq x_{1}\geq 0\}.$$
It is easy to see that $K=K_{\infty}$,	$f^{\infty}_{1}(x_{1}, x_{2})=x^{2}_{1}$, and $f^{\infty}_{2}(x_{1}, x_{2})=x^{3}_{2}$. Clearly, $SOL^{s}(K_{\infty}, f^{\infty})=\{0\}$.
On the other hand, $SOL^{s}(K, f)$ is unbounded since
$$\{(x_{1}, x_{2})\in K: x_{1}=0, x_{2}\geq 0\}\subseteq SOL^{s}(K, f).$$
	\end{example}

The following result gives a Frank-Wolf type theorem for  $text{ PVOP } (K, f)$ under the weak regularity condition.

\begin{corollary}\label{wexist}
If $\mathrm{ PVOP } (K, f)$ is weakly regular and  $f: \mathbf{R}^{n}\mapsto \mathbf{R^s}$ is bounded from below on $K$, then $SOL^{s}(K, f)$ is nonempty.
\end{corollary}
{\it Proof}
It  follows directly from Proposition \ref{weakpro}, Remark \ref{RE} and Theorem \ref{exist0}.
\qed

The following result shows that the existence of the Pareto efficient solutions of a polynomial vector optimization problem  is preserved  when  the vector objective function is perturbed by a lower degree polynomial.

\begin{theorem}\label{existper}
Assume that one of the following conditions hold:
\begin{itemize}
\item[(i)] $SOL^{s}(K_{\infty}, f^{\infty})=\{0\}$.
\item[(ii)] there exists $i_{0}\in \{1, \dots, s\}$ such that $SOL(K_{\infty}, f^{\infty}_{i_{0}})=\{0\}$.
\end{itemize}
Then $SOL^{s}(K, f+g)$ is nonempty for all $g=(g_{1},  \dots, g_{s})\in \mathbf{P}_{\mathbf{d}}$ with  $\mathrm{ deg }\ g_{i}<\mathrm{ deg }\ f_{i}$, $i=1,  \dots, s$.
\end{theorem}
{\it Proof}
It follows directly from  Proposition \ref{perturb}\emph{(i)} and Theorem \ref{exist0}.
 \qed

\subsection{The strong regularity case}

In this subsection, we study nonemptiness and boundedness of the solution sets of $\text{ PVOP }(K, f)$ under the strong regularity condition. Next, we give a necessary condition for the existence of the weak Pareto efficient solutions of  $\text{ PVOP }(K, f)$.
\begin{theorem}\label{special1}
If $SOL^{w}(K_{\infty}, f^{\infty})=\emptyset$, then $SOL^{w}(K, f)=\emptyset$.
\end{theorem}
{\it Proof}
Suppose that there exists $x_{0}\in SOL^{w}(K, f)$.  By Proposition \ref{necessary}, $0\notin SOL^{w}(K_{\infty}, f^{\infty})$. Then there exists $v_{1}\in K_{\infty}\backslash\{0\}$ such that $$f^{\infty}(v_{1})-f^{\infty}(0)=f^{\infty}(v_{1})\in -\mbox{ int }\mathbf{R}^{s}_{+}.$$
This means that
\begin{equation}\label{fypa2}
f^{\infty}_{i}(v_{1})<f^{\infty}_{i}(0)=0,\quad i=1,  \dots, s.
\end{equation}
Since  $v_{1}\in K_{\infty}$, there exist $t_k >0$ with $t_{k}\rightarrow +\infty$ and $x_{k}\in K$ such that $t^{-1}_{k}x_{k}\rightarrow v_{1}$ as $k\rightarrow +\infty$.
Since $x_{0}\in SOL^{w}(K, f)$ and $x_k\in K$,  there exists $i_{x_k}\in \{1, \dots, s\}$ such that
$$
f_{i_{x_k}}(x_k)-f_{i_{x_k}}(x_{0})\geq 0, \quad \forall k.
$$
Without loss of generality, we can assume that   there exists $i_{0}\in \{1, \dots, s\}$ such that
$$
	f_{i_{0}}(x_{k})-f_{i_{0}}(x_{0})\geq 0.
$$
Dividing the both sides of the above inequality by $t_{k}^{d_{i_{0}}}$ and letting $k\rightarrow +\infty$, we have $$f^{\infty}_{i_{0}}(v_{1})\geq 0,$$
 a contradiction to (\ref{fypa2}).
\qed

\begin{remark}
By Theorem \ref{special1}, $SOL^{w}(K, f)\neq\emptyset$ implies that  $SOL^{w}(K_{\infty}, f^{\infty})\neq\emptyset$. The following example shows that the converse does not hold in general.
\end{remark}

\begin{example}
Consider the vector polynomial  $f=(f_{1}, f_{2})$ with
$$f_{1}(x_{1}, x_{2})=(x^{4}_{1}x^{4}_{2}-1)^{2}+x^{4}_{1}, f_{2}(x_{1}, x_{2})=(x^{2}_{1}x^{2}_{2}-1)^{2}+x^{2}_{1}$$
and $K=\mathbf{R}^{n}$. It is easy to verify that $0\in SOL^{w}(K_{\infty}, f^{\infty})$. On the other hand, $f_{1}>0$ and $f_{2}>0$, but $f(\frac{1}{n}, n)=(\frac{1}{n^{4}}, \frac{1}{n^{2}})\to (0, 0)$ as $n\to +\infty$. This implies $SOL^{w}(K, f)=\emptyset$.
\end{example}

As a direct consequence of  Proposition \ref{necessary} and Theorem \ref{special1}, we have the following result.

\begin{theorem}\label{th-fyp1}
If  $0\notin SOL^{w}(K_{\infty}, f^{\infty})$, then  $SOL^{w}(K, f)=\emptyset$.
\end{theorem}

\begin{remark}
By Theorem \ref{special1} (or Theorem \ref{th-fyp1}), 	$SOL^{w}(K, f)\ne\emptyset$ implies  $SOL^{w}(K_{\infty}, f^{\infty})\ne\emptyset$. The following example shows that  $SOL^{s}(K_{\infty}, f^{\infty})$ may be empty when 	$SOL^{s}(K, f)\ne\emptyset$.
\end{remark}

\begin{example}
		Consider the vector polynomial $f=(f_{1}, f_{2})$ with
		$$f_{1}(x_{1}, x_{2})=-x^{2}_{2}-1, f_{2}(x_{1}, x_{2})=-x^{3}_{1}+x_{2}+1$$
		and
$$K=\{(x_{1}, x_{2})\in \mathbf{R}^{2}: 0\leq x_{1}\leq 1, x_{2}\geq 2\}.$$	
Then $K_{\infty}=\{(x_{1}, x_{2})\in \mathbf{R}^{2}: x_{1}=0, x_{2}\geq 0\}$ and $f^{\infty}(x_{1}, x_{2})=(f^{\infty}_{1}(x_{1}, x_{2}), f^{\infty}_{2}(x_{1}, x_{2}))=(-x^{2}_{2}, -x^{3}_{1})$. It is easy to verify that  $SOL^{s}(K_{\infty}, f^{\infty})=\emptyset$ and  $(1, 2)\in SOL^{s}(K, f)$.
	\end{example}

The following result shows that $SOL^{w}(K_{\infty}, f^{\infty})=\{0\}$ is  sufficient for the existence of the Pareto efficient solutions as well as boundedness of the weak Pareto efficient solution set.

\begin{theorem}\label{existt}
 If $SOL^{w}(K_{\infty}, f^{\infty})=\{0\}$, then $SOL^{s}(K, f)$ is nonempty and $SOL^{w}(K, f)$ is compact.
\end{theorem}
{\it Proof}
$SOL^{s}(K, f)\ne\emptyset$ follows from  Theorem \ref{gpro}\emph{(i)} and Theorem \ref{exist0}\emph{(ii)}. The closedness of $SOL^{w}(K, f)$ is clear. Next we prove that $SOL^{w}(K, f)$ is bounded. If not,  then there exists $y_{k}\in SOL^{w}(K, f)$ such that $\|y_{k}\|\rightarrow +\infty$ as $k\rightarrow +\infty$. Without loss of generality, assume that $\|y_{k}\|\neq 0$ and $\frac{y_{k}}{\|y_{k}\|}\rightarrow v\in K_{\infty}\backslash\{0\}$. Let $y\in K$. Since $y_{k}\in SOL^{w}(K, f)$, there exists $i_{y_{k}}\in \{1,  \dots, s\}$ such that
	$$f_{i_{y_{k}}}(y)-f_{i_{y_{k}}}(y_{k})\geq 0.$$
Since the set $\{1,  \dots, s\}$ is finite, without loss of generality, we can assume that there exists $i_{0}\in \{1,  \dots, s\}$ such that
	$$f_{i_0}(y)-f_{i_0}(y_{k})\geq 0,\quad \forall k.$$
Dividing the both sides of the above inequality by $\|y_{k}\|^{d_{i_{0}}}$ and letting $k\to +\infty$, we have $f^{\infty}_{i_{0}}(v)\leq f^{\infty}_{i_{0}}(0)=0$. This implies $v\in SOL(K_{\infty}, f^{\infty}_{i_{0}})\backslash\{0\}$, a contradiction.
\qed

As a consequence of Theorem \ref{existt} and Theorem \ref{gpro}\emph{(i)}, we have the following result.

\begin{corollary}\label{simple}
 	If $SOL(K_{\infty}, f^{\infty}_{i})= \{0\}$ for each $i\in \{1, 2, \dots, s\}$, then $SOL^{s}(K, f)$ is nonempty and $SOL^{w}(K, f)$ is bounded.
 \end{corollary}

Now we give an example to illustrate the conclusion of Corollary \ref{simple}.

\begin{example}\label{1}
	Consider the vector polynomial $f=(f_{1}, f_{2})$ with
 $$f_{1}(x_{1}, x_{2})=x^{3}_{2}-x^{2}_{1}-x_{1}x_{2}+1, f_{2}(x_{1}, x_{2})=x^{2}_{2}-x_{1}-1$$ and
  $$K=\{(x_{1}, x_{2})\in \mathbf{R}^{2}: x_{1}-1\geq 0, x_{2}-x_{1}+1\geq 0,e^{x_{1}-1}-x_{2}\geq 0\}.$$
It is easy to verify that $f^{\infty}_{1}(x_{1}, x_{2})=x^{3}_{2},f^{\infty}_{2}(x_{1}, x_{2})=x^{2}_{2}$,
$$K_{\infty}=\{(x_{1}, x_{2})\in \mathbf{R}^{2}: x_{1}\geq 0, x_{2}-x_{1}\geq 0\},$$
$SOL(K_{\infty}, f^{\infty}_{1})=\{(0, 0)\}$, and $SOL(K_{\infty}, f^{\infty}_{2})=\{(0, 0)\}$. By Corollary \ref{simple},  $SOL^{s}(K, f)$ is nonempty and $SOL^{w}(K, f)$ is compact.  It is worth mentioning that \cite[Theorem 4.1]{DTN} and \cite[ Theorem 3.1]{LGJ} cannot be applied in this example since $f$ is non-convex on $\mathbf{R}^{2}$ and $K$ is neither convex nor semi-algebraic set.
\end{example}

\begin{remark}\label{rmk6.3}
By Theorem \ref{exist0}\emph{(ii)}, $SOL^{s}(K, f)$ is nonempty when there exists some $i_0\in \{1,\cdots, s\}$ such that $SOL(K_{\infty}, f^{\infty}_{i_0})=\{0\}$. If $SOL(K_{\infty}, f^{\infty}_{i})=\{0\}$ for all  $i\in \{1, \cdots, s\}$, then  $SOL^{s}(K, f)$ is nonempty and bounded(by Corollary \ref{simple}). The following example shows that $SOL^{s}(K, f)$ may  be unbounded when there exists $\{i_0, j_0\}\subset \{1,\cdots, s\}$ such that $SOL(K_{\infty}, f^{\infty}_{i_{0}})=\{0\}$ and $SOL(K_{\infty}, f^{\infty}_{j_{0}})\ne\{0\}$.
\end{remark}

\begin{example}
		Consider the vector polynomial $f=(f_{1}, f_{2})$ with
$$f_{1}(x_{1}, x_{2})=x^{3}_{1}, \quad f_{2}(x_{1}, x_{2})=-x^{2}_{1}+x_{2}$$ and
$$K=\{(x_{1}, x_{2})\in \mathbf{R}^{2}: x_{2}\geq 0, x_{1}-x_{2}\geq 0\}.$$
It is easy to verify that $K=K_{\infty}$, $f^{\infty}_{1}(x_{1}, x_{2})=x^{3}_{1},f^{\infty}_{2}(x_{1}, x_{2})=-x^{2}_{1}$, $SOL(K_{\infty}, f^{\infty}_{1})=\{0\}$, $SOL(K_{\infty}, f^{\infty}_{2})=\emptyset$.
On the other hand,  $SOL^{s}(K, f)$ is unbounded since
$$\{(x_{1}, x_{2})\in K: x_{1}\geq 0, x_{2}=0\}\subseteq SOL^{s}(K, f).$$
\end{example}

It has been shown in Corollary \ref{wexist} that  $\text{ PVOP } (K,f)$ with  weak regularity admits a weak Pareto efficient solution provided that $f$ is bounded from below on $K$.  In the following Frank-Wolf type theorem, we further prove the existence of Pareto efficient solutions and compactness of  the weak Pareto efficient solution set if we strengthen weak regularity by  strong regularity.

\begin{corollary}\label{exist}
If \ $\mathrm{ PVOP }(K, f)$ is  strongly regular and  $f: \mathbf{R}^{n}\mapsto \mathbf{R^s}$ is bounded from below on $K$, then $SOL^{s}(K, f)$ is nonempty and $SOL^{w}(K, f)$ is compact.
\end{corollary}
{\it Proof}
 It follows directly from Proposition \ref{eqcondition} and Theorem \ref{existt}.
\qed

\begin{remark}
Corollary \ref{exist} extends \cite[Theorem 3.1]{HV1} to the vector case.
\end{remark}

The following example shows that the converse of Theorem \ref{existt}, Corollary \ref{simple} and Corollary \ref{exist} does not hold in general.

\begin{example}
	Consider the polynomial $f=(f_{1}, f_{2})$ with
$$f_{1}(x_{1}, x_{2})=x_{1}x_{2}+1,\quad f_{2}(x_{1}, x_{2})=x_{1}x_{2}+x_{1}-1$$ and
$$K=\{(x_{1}, x_{2})\in \mathbf{R}^{2}: x_{1}\geq 1, x_{2}\geq 1\}.$$
	Then $f^{\infty}_{1}(x_{1}, x_{2})=f^{\infty}_{2}(x_{1}, x_{2})=x_{1}x_{2}$ and  $K_{\infty}=\mathbf{R}^{2}_{+}$. It is easy to verify that $SOL^{s}(K, f)=SOL^{w}(K, f)=\{(1, 1)\}$. However, $SOL(K_{\infty}, f^{\infty}_{1})=SOL(K_{\infty}, f^{\infty}_{2})=\{(x_{1}, x_{2})\in \mathbf{R}^{2}_{+}: x_{1}x_{2}=0\}$ are unbounded.
\end{example}

The following result shows that some properties of  the  solution sets of a strongly regular polynomial vector optimization problem  are preserved when its objective polynomial  is perturbed by  a lower degree polynomial.

\begin{theorem}\label{empper} Let $g=(g_1,\cdots, g_s): \mathbf{R}^{n} \mapsto \mathbf{R}^{s}$ be a vector polynomial with $\mathrm{ deg }$\ $g_i<\mathrm{ deg }$\ $ f_i, i=1,\cdots, s$. Then the following conclusions hold:
\begin{itemize}
\item[(i)] If $SOL^{w}(K_{\infty}, f^{\infty})=\emptyset$, then $SOL^{w}(K, f+g)=\emptyset$.
\item[(ii)] If $SOL^{w}(K_{\infty}, f^{\infty})=\{0\}$, then $SOL^{s}(K, f+g)$ is nonempty and $SOL^{w}(K, f+g)$ is compact.
\end{itemize}
\end{theorem}
{\it Proof}
\emph{(i)} follows from  Proposition \ref{perturb}\emph{(iv)} and Theorem \ref{special1}, and \emph{(ii)} follows from Proposition \ref{perturb}\emph{(ii)} and Theorem \ref{existt}.
\qed

\section{Stability Analysis}

In this section we investigate the solution stability of a regular polynomial vector optimization problem. First, we shows that some properties of the solution sets of a strongly regular polynomial vector optimization problem are stable under a small perturbation.

\begin{theorem}\label{localper} The following conclusions hold:
\begin{itemize}
\item[(i)] If $SOL^{w}(K_{\infty}, f^{\infty})=\emptyset$, then there exists $\epsilon>0$ such that $SOL^{w}(K, f+g)=\emptyset$ for all $g\in \mathbf{P}_{\mathbf{d}}$ satisfying $\|g\|<\epsilon$.
\item[(ii)] If $SOL^{w}(K_{\infty}, f^{\infty})=\{0\}$, then there exists $\epsilon>0$ such that $SOL^{s}(K, f+g)$ is nonempty and $SOL^{w}(K, f+g)$ is compact for all $g\in \mathbf{P}_{\mathbf{d}}$ satisfying $\|g\|<\epsilon$.
\end{itemize}
\end{theorem}
{\it Proof}
\emph{(i)} follows from Theorem \ref{local regu} and Theorem \ref{special1}	and \emph{(ii)} follows from Theorem \ref{local regu} and Theorem \ref{existt}.
\qed

In the sequel we shall investigate the local boundedness and upper semicontinuity of the weak Pareto efficient solution mapping.  Recall that a set-valued mapping $T:\mathbf{R}^{\kappa}\rightrightarrows 2^{\mathbf{R}^{n}}$ is said to be upper semi-continuous at $x$ if for any neighborhood $U$ of $T(x)$, there exists a neighborhood $V$ of $x$ such that $T(y)\subseteq U$ for any $y\in V$. A set-valued mapping $T: \mathbf{R}^{\kappa}\rightrightarrows 2^{\mathbf{R}^{n}}$ is locally bounded at $x$ if there exists an open neighborhood $V$ of $x$ such that $\cup_{y\in V}T(y)$ is bounded.

\begin{lemma}\label{lems6} \cite[Theorem 5.7 and Theorem 5.19]{RC}
If $T:\mathbf{R}^{\kappa}\rightrightarrows 2^{\mathbf{R}^{n}}$ is locally bounded at $x$  and  $T$ has a closed  graph $\mathrm{Gph}(T)$, where $\mathrm{Gph}(T)=\{(x, y)\in \mathbf{R}^{\kappa}\times\mathbf{R}^{n}: y\in T(x)\}$, then  $T$ is upper semi-continuous at $x$.
\end{lemma}

\begin{theorem}\label{uppercon}
	Assume that $K$ is convex and $\mathrm{PVOP }(K, f)$ is strongly  regular. Then the
following conclusions hold:
\begin{itemize}
\item[(i)] $SOL^{w}(K, \cdot)$ is locally bounded at $f$, i.e., there exists $\delta>0$ such that $$Q_{\delta}=\cup_{h\in \mathbf{B}(0, \delta)}SOL^{w}(K, f+h)$$
is bounded, where the $\mathbf{B}(0, \delta)$ is an open ball in $\mathbf{P}_{\mathbf{d}}$ with center at $0$ and radius $\delta>0$.
\item[(ii)] $SOL^{w}(K, \cdot)$ is upper semi-continuous on $\mathbf{GR}^{\mathbf{d}}_{s}$.
\end{itemize}
\end{theorem}
{\it Proof}
	\emph{(i)}  Since   $\mathbf{GR}^{\mathbf{d}}_{s}$ is open in $\mathbf{P}_{\mathbf{d}}$ (by Proposition \ref{open}) and $f\in\mathbf{GR}^{\mathbf{d}}_{s}$, there exists $\delta>0$ such that $f+\bar{\mathbf{B}}(0, \delta)\subseteq\mathbf{GR}^{\mathbf{d}}_{s}$, where $\bar{\mathbf{B}}(0, \delta)$ denotes the closure of  $\mathbf{B}(0, \delta)$.  Suppose on the contrary that $Q_{\delta}=SOL^{w}(K, f+\mathbf{B}(0, \delta))$ is unbounded. Then there exist $\{h^{k}\}\subseteq\mathbf{B}(0, \delta)$ and $x_{k}\in SOL^{w}(K, f+h^{k})$ such that $\|x_k\|\to +\infty$ as $k\to +\infty$.  Without loss of generality, we may assume that  $\frac{x_{k}}{\|x_{k}\|}\to x^{*}\in K_{\infty}\backslash \{0\}$ and $h^{k}\to h\in \bar{\mathbf{B}}(0, \delta)$.
Let $x\in K$ and $v\in K_{\infty}$. Then $x+\|x_{k}\|v\in K$ for all $k$. Since  $x_{k}\in SOL^{w}(K, f+h^{k})$, we have
	$$(f+h^{k})(x+\|x_{k}\|v)-(f+h^{k})(x_{k})\notin -\mbox{ int }\mathbf{R}^{s}_{+}.$$
Then for each $k$, there exists $i_{k, v}\in \{1, \dots, s\}$ such that
	$$(f_{i_{k, v}}+h^{k}_{i_{k, v}})(x+\|x_{k}\|v)-(f_{i_{k, v}}+h^{k}_{i_{k, v}})(x_{k})\geq 0.$$
Since $\{1,  \dots, s\}$ is finite, we can assume that there exists $i_{0, v}\in \{1, \dots, s\}$ such that
$$(f_{i_{0, v}}+h^{k}_{i_{0, v}})(x+\|x_{k}\|v)-(f_{i_{0, v}}+h^{k}_{i_{0, v}})(x_{k})\geq 0,\quad\forall k.$$
Dividing the both sides of the above inequality by $\|x_{k}\|^{d_{i_{0, v}}}$ and letting $k\to +\infty$, we get
$$(f_{i_{0, v }}+h_{i_{0, v}})^{\infty}(v)\geq (f_{i_{0, v}}+h_{i_{0, v}})^{\infty}(x^{*}).$$
Since  $v\in K_{\infty}$ is arbitrary, we have $x^{*}\in SOL^{w}(K_{\infty}, (f+h)^{\infty})\backslash\{0\}$. On the other hand, since $f+h\in f+\bar{\mathbf{B}}(0, \delta)\subseteq \mathbf{GR}^{\mathbf{d}}_{s}$,  by Proposition \ref{cones} we have $x^{*}\notin SOL^{w}(K_{\infty}, (f+h)^{\infty})$, a contradiction.  Thus, $SOL^{w}(K, \cdot)$ is locally bounded at $f$.
	
    \emph{(ii)} Let $f\in \mathbf{GR}^{\mathbf{d}}_{s}$. By \emph{(i)} and Lemma \ref{lems6}, we only need to prove that the graph  $\mbox{Gph}(SOL^{w}(K, \cdot))$ of $SOL^{w}(K, \cdot)$  is closed in $\mathbf{P}_{\mathbf{d}}\times\mathbf{R}^{n}$. Let $(f^{k}, y_{k})\to (f, y)$ as $k\to +\infty$ with $ y_{k}\in SOL^{w}(K, f^{k})$. Then for any $z\in K$, we have
$$f^{k}(z)-f^{k}(y_{k})\notin -\mbox{ int }\mathbf{R}^{s}_{+}.$$
Letting $k\to +\infty$, we have	
$$f(z)-f(y)\not\in -\mbox{ int }\mathbf{R}^{s}_{+},\quad \forall z\in K,$$
which means  $y\in SOL^{w}(K, f)$.  Thus, $\mbox{Gph}(SOL^{w}(K,\cdot))$ is closed in $\mathbf{P}_{\mathbf{d}}\times\mathbf{R}^{n}$..
\qed

\section{Conclusion}
In this paper we extend the concept of regularity due to Hieu \cite{HV1} to  the polynomial vector optimization problem. Under  regularity conditions, we investigate nonemptiness and boundedness of the solution sets of a non-convex polynomial vector optimization problem on a nonempty closed set (not necessarily  semi-algebraic set).  As a consequence, we derive two Frank-Wolfe type theorems for a non-convex polynomial vector optimization problem. We also discuss the solution stability.  Our results extend and improve the corresponding results of \cite{DTN,BT1,BT2,HV1,LGJ}.



\begin{thebibliography}{1}
\bibitem{HHV} Ha, H.V., Pham, T.S.(eds.): Genericity in Polynomial Optimization. World Science Publishing, Singapore (2017)

\bibitem{BR} Benedetti, R., Risler, J.J.(eds.): Real Algebraic and Semi-Algebraic Sets. Hermann, Paris (1990)

\bibitem{AQ} Ansari, Q.H., Lalitha, C.S., Mehta, M.(eds.): Generalized Convexity, Nonsmooth Variational Inequalities, and Nonsmooth Optimization. Chapman and Hall/CRC, New York (2013)

\bibitem{HV1} Hieu, V.T.: A regularity condition in polynomial optimization. arXiv:1808.06100

\bibitem{HV2} Hieu, V.T., Wei, Y.M., Yao, J.C.: Notes on the optimization problems corresponding to polynomial complementarity problems. J. Optim. Theory Appl. \textbf{184}, 687--695 (2020)

\bibitem{DTN} Kim, D.S., Pham, T.S., Tuyen, V.T.: On the existence of Pareto solutions for polynomial vector optimization problems. Math. Program. \textbf{177}, 321--341 (2019)

\bibitem{LGJ} Jiao, L.G., Lee, J.H., Sisarat, N.: Multi-objective convex polynomial optimization and semidefinite programming relaxations. arXiv:1903.10137

\bibitem{RC} Rockafellar, R.T., Wets, R.J-B.(eds.): Variational Analysis. Springer, Berlin (2009)

\bibitem{AA} Auslender, A., Teboulle, M.(eds.): Asymptotic Cones and Functions in Optimization and Variational Inequalities. Springer, New York (2003)

\bibitem{BT1} Bao, B.T., Mordukhovich, B.S.: Variational principles for set-valued mappings with applications to multiobjective optimization. Control Cybern. \textbf{36}, 531--562 (2007)

\bibitem{BT2} Bao, B.T., Mordukhovich, B.S.: Relative Pareto minimizers for multiobjective problems: existence and optimality conditions. Math. Program. \textbf{122}, 301--347 (2010)

\bibitem{MG} Giannessi, F., Mastroeni, G., Pellegrini, L.: On the theory of vector optimization and variational inequalities. Image space analysis and separation. In: Giannessi, F. (ed.): Vector Variational Inequalities and Vector Equilibria. Nonconvex Optimization and Its Applications, vol. 38, Springer, Boston, MA. (2000)

\bibitem{FYP} Huang, N.J. and Fang, Y.P.: The upper semicontinuity of the solution maps in vector implicit quasicomplementarity problems of type $R_0$. Appl. Math. Lett. \textbf{16}, 1151--1156 (2003)

\bibitem{FangH}Fang, Y.P. and Huang, N.J.: On the upper semi-continuity of the solution map to the vertical implicit homogeneous complementarity problem of type $R_0$. Positivity \textbf{10}, 95--104 (2006)

\bibitem{Oett}Oettli, W., Yen, N.D.: Quasicomplementarity problems of type $R_0$, J. Optim. Theory Appl. \textbf{89}, 467--474 (1996)

\bibitem{Lop} R. L\'{o}pez: Stability results for polyhedral complementarity problems. Comput. Math. Appl.  \textbf{58}. 1475--1486 (2009)

\bibitem{Gowda} Gowda, M. S., Sossa, D.: Weakly homogeneous variational inequalities and
solvability of nonlinear equations over cones. Math. Program., Ser. A, \textbf{177}, 149--171 (2019)

\bibitem{Yangxq}Yang, X. Q.: Vector complementarity and minimal element problems.
J. Optim. Theory Appl. \textbf{77}, 483--495 (1993)

\bibitem{CGY} Chen, G.Y., Yang, X.Q.: The vector complementary problem and its equivalence with weak minimal elements in ordered spaces. J. Math. Anal. Appl. \textbf{153}, 136--158 (1990)

\bibitem{LKY} Lee, G.M., Kim, D.S., Lee, B.S., Yen N.D.: Vector variational inequality as a tool for studying vector optimization problems. In: Giannessi, F.(ed.) Vector Variational Inequalities and Vector Equilibria. Nonconvex Optimization and Its Applications, vol. 38, Springer, Boston, MA. (2000)


\end{thebibliography}
\end{document}